\documentclass[runningheads]{llncs}
\usepackage{amssymb,amstext,amsmath,mathabx}
\usepackage[utf8]{inputenc}
\usepackage{hyperref}
\usepackage{verbatim}
\usepackage{mathtools}
\usepackage{caption,color}
\newcommand*{\normally}{\mathrel{\ooalign{$|$\hfil\cr\kern+1pt$\thicksim$}}} 
\newcommand*{\nnormally}{\mathrel{\ooalign{$|$\hfil\cr\kern+1pt$\thicksim$}\negthickspace \negthickspace /} }

\def\G{\mathcal{G}}
\def\H{\mathcal{H}}
\def\F{\mathcal{F}}

\def\F{\mathcal{F}}
\def\P{\mathcal{P}}

\def\prev{\mathbb{P}}

\newcommand{\no}[1]{\widebar{#1}}
\newcommand*{\lsim}{\mathord{\thicksim}}
\renewcommand{\bot}{\emptyset}
\title{Interpreting connexive principles in coherence-based probability logic}
\author{Niki Pfeifer\inst{1}\thanks{Both authors contributed equally to the article and are listed alphabetically.
}\,\thanks{Supported by the BMBF project  01UL1906X.} 
\and  Giuseppe Sanfilippo\inst{2}$^\star$\;\thanks{Member of the GNAMPA Research Group 
and partially supported by the INdAM–GNAMPA Project 2020 Grant U-UFMBAZ-2020-000819.}
}

\institute{Department of Philosophy, University of Regensburg, Germany
   \\ \email{niki.pfeifer@ur.de}
  \and
  Department of Mathematics and Computer Science,
University of Palermo, Italy  
 \\ \email{giuseppe.sanfilippo@unipa.it}
}
\authorrunning{Niki Pfeifer and Giuseppe Sanfilippo}
\begin{document}
\maketitle
\begin{abstract}
We present probabilistic approaches to check the validity of selected connexive principles within the setting of coherence.  Connexive logics emerged from the intuition that conditionals of the
form \emph{If $\lsim A$, then $A$}, should not hold, since the conditional's antecedent $\lsim A$ contradicts its consequent $A$. Our approach covers this intuition by observing  that for an event A the only coherent probability assessment on the conditional event $A|\no{A}$  is $p(A|\no{A})=0$. Moreover, connexive logics aim to capture the intuition that conditionals should express some ``connection"
between the antecedent and the consequent or, in terms of inferences, validity should require some connection between the premise set and the conclusion. This intuition is covered by a number of principles, a selection of which we analyze in our contribution. We present two approaches to connexivity within coherence-based probability logic. Specifically, we analyze connections between antecedents and consequents firstly, in terms of probabilistic constraints on conditional events (in the sense of defaults, or negated defaults) and secondly,  in terms of constraints on compounds of conditionals and iterated conditionals. After developing different notions of negations and notions of validity, we analyze the following connexive principles within both approaches: Aristotle's Theses, Aristotle's Second Thesis, Abelard's First Principle and selected versions of Boethius' Theses. We conclude by remarking that coherence-based probability logic offers a rich language to investigate the validity of various connexive principles.
\end{abstract}
\section{Introduction}
We present  probabilistic approaches to check the validity of selected connexive principles  within the setting of coherence.  Connexive logics emerged from the intuition that conditionals of the
form \emph{if not-$A$, then $A$}, denoted by  $\sim A \rightarrow  A$, should not hold, since the
conditional's antecedent \emph{not-$A$} contradicts its consequent $A$. Indeed, experimental psychological data show that people believe that sentences of the form \emph{if not-$A$, then $A$} are false (e.g., \cite{pfeifer12x,pfeiferPLip}), which supports the psychological plausibility of this intuition.   Connexive principles were
developed to rule out such self-contradictory conditionals  (for overviews, see
e.g.,\cite{mccall12,wansing20}). 
Many of
these principles can be traced back to antiquity or the middle ages,
which is reflected by the  names of these principles, for example, Aristotle's Thesis or Abelard's First Principle (see Table \ref{TAB:CP}).
\begin{table}[!h]	
	\centering \vspace{-1em}
		\begin{tabular}{lll}	
			\hline 
			Name & Abbreviation
			&	
			Connexive principle   \\\hline	&\\[-1em] 
			Aristotle's Thesis & (AT) &$\lsim (\lsim A \rightarrow A)$ \\	
			
			&\\[-1em] 
			
			Aristotle's Thesis$'$&(AT$'$) &$\lsim ( A \rightarrow \lsim A)$ \\
			
			&\\[-1em] 
			
			Abelard’s First Principle&(AB) &$\lsim((A \rightarrow B) \wedge (A \rightarrow\lsim B))$\\
			
			&\\[-1em] 
			
			Aristotle's Second Thesis&(AS) & $\lsim((A \rightarrow B) \wedge (\lsim A \rightarrow B))$\\
			
			&\\[-1em] 
			
			Boethius' Thesis&(BT) & $(A \rightarrow B) \rightarrow \lsim (A
			\rightarrow \lsim B)$\\
			
			&\\[-1em] 
			
			Boethius' Thesis$'$&(BT$'$) & $(A \rightarrow \lsim B) \rightarrow
			\lsim(A \rightarrow B)$ \\
			
			&\\[-1em] 
			
			Reversed Boethius' Thesis&(RBT) & $\lsim (A
			\rightarrow \lsim B) \rightarrow(A \rightarrow B)$ \\

			&\\[-1em] 
			
			Reversed Boethius' Thesis$'$&(RBT$'$) & $ \lsim(A \rightarrow B)\rightarrow(A \rightarrow \lsim B)$\\
			
			&\\[-1em] 
			
			Boethius Variation 3&(B3) & $(A \rightarrow B) \rightarrow  \lsim (\lsim A
			\rightarrow  B)$ \\
			
			&\\[-1em] 
			
			Boethius Variation 4&(B4)& $(\sim A \rightarrow B) \rightarrow  \lsim (A
			\rightarrow  B)$\\
			
			&\\[-1.1em] 
			
			\hline 
		\end{tabular}
	\caption{\label{TAB:CP}Selected connexive principles (see also \cite{wansing20}).}
\end{table} \vspace{-2.3em}
In classical logic, however, 
Aristotle's Thesis, i.e. $\sim (\sim A \rightarrow  A)$, 
is not a theorem  since the corresponding material conditional  is contingent because $\sim ( \sim \sim A \vee A)$ is logically equivalent to $\sim A$ (which is not necessarily true).  
Moreover, connexive logics
aim to capture the intuition  that conditionals should express some ``connection''
between the antecedent and the consequent or, in terms of inferences,
validity should require some connection between the premise set and
the conclusion.

The connexive intuition that conditionals of the
form \emph{if not-$A$, then $A$} should not hold  is covered in  subjective probability theory.
Specifically, we cover this intuition by the observation that for any  event $A$, with $\no{A}\neq \bot$, the  only coherent  assessment on the conditional event  $A|\no{A}$  is $p(A|\no{A})=0$.

The aim of our contribution is to investigate selected connexive
principles within the framework of \emph{coherence-based probability
  logic}. The coherence
approach to (subjective) probability  was
originated by Bruno de Finetti (see, e.g., \cite{definetti31,definetti74}) and  has been generalised  to the
conditional probability 
and to  previsions of conditional random quantities
(see, e.g., \cite{BeRR98,biazzo00,CaLS07,coletti02,gilio13,GiSa19,Holz85,lad96,Rega85,WaPV04}). 
In the present framework, we present two approaches to connexivity within coherence-based probability logic. In the first approach  we analyze connections between antecedents and consequents in terms of probabilistic constraints on conditional events (in the sense of defaults or negated defaults \cite{gilio16,PS17SH,PfSa18,PfSa19}). In the second approach, based the
recently developed more general framework of compounds of conditionals and iterated conditionals (\cite{GiSa13c,GiSa13a,GiSa14,GiSa19,GiSa21}), we define these connections in terms of constraints on suitable conditional random quantities.
After developing different notions of negations and notions of validity, we analyze the connexive principles given in Table \ref{TAB:CP} within both approaches.

The  coherence principle  plays a key role  in probabilistic reasoning and allows for probabilistic inferences of a further conditional event (the  
conclusion) from   any  coherent probabilistic assessment on an arbitrary family of
conditional events (the premises).
Moreover,  coherence  is a more general approach to conditional probabilities compared to approaches which requires positive probability for the conditioning events.  In  standard  approaches to probability the  conditional probability $p(C|A)$ is defined by the ratio $p(A\wedge C)/p(A)$, which requires  positive probability of the conditioning event,  $p(A)>0$. 
 However,   in the framework of coherence, conditional probability
 $p(C|A)$, as a degree of belief,  is a primitive notion and it is properly defined  even if
 the  conditioning event  has  probability zero, i.e.,
 $p(A)=0$. Analogously, within coherence, previsions of conditional random quantities, are primitive and  properly defined even if the conditioning event has probability zero.
 Therefore,  coherence  is a more general approach to
 conditional probabilities compared to approaches which requires
 positive probability for the conditioning events. The only
 requirement is that the conditioning event must be logically
 possible. Thus, although $p(C|A)$ is well defined even if
 $p(A)=0$, it is undefined if $A\equiv \bot$ (where $\bot$ denotes a
 logical contradiction). This is in line with the reading that Boethius
 and Aristotle thought  that  principles like (BT) and (AT), respectively, hold
 only when the conditional's antecedent is
 possible  (see \cite{lenzen20} who argues that the ``ancient logicians most
 likely meant their theses as applicable only to `normal' conditionals
 with antecedents which are not self-contradictory'';
 p. 16). It is also in line with the Ramsey test, which is
 expressed in his famous footnote: 
``If two people are arguing `If $A$ will $C$?' and are both in doubt
as to $A$, they are adding $A$ hypothetically to their stock of
knowledge and arguing on that basis about $C$; so that in a sense `If
$A$, $C$' and `If $A$, $\no{C}$' are contradictories. We can say
they are fixing their degrees of belief in $C$ given $A$. If $A$ turns
out false, these degrees of belief are rendered \emph{void}'' \cite[p. 155,
we adjusted the notation]{ramsey29}.
 The quantitative interpretation of the Ramsey test became a cornerstone of the
 conditional probability interpretation of conditionals. Adding a contradiction to your stock of knowledge does not
make sense (as, traditionally, knowledge implies truth). Moreover, Ramsey's thought that conditionals with contradicting
consequents $C$ and $\no{C}$ contradict each other coincides with the
underlying intuition of  (AB).
\section{Preliminary Notions and Results}
Given two events $A$ and  $H$, with $H\neq \bot$ (where $\bot$ denotes the impossible event),  the \emph{conditional event} $A|H$ (read: $A$ given $H$) is defined as a three-valued logical entity
which is \emph{true} if $AH$ (i.e., $A\wedge H$) is true, \emph{false} if $\no{A}H$ is true, and  \emph{void} if $H$ is false.  We observe that $A|H$  assumes the logical value (\emph{true} or \emph{false}) of $A$, when $H$ is true, and it is \emph{void}, otherwise.
There is a long history of how to deal with negations (see, e.g., \cite{horn01}).  In our context, the  negation  of the conditional event ``$A$ given $H$'', denoted by $\overline{A|H}$, is the conditional event 
$\no{A}|H$, that is   
``the negation of  $A$'' given $H$.  We use the inner negation to preserve for conditional events the usual property of negating unconditional events: $p(\no{A})=1-p(A)$.
In the subjective approach to probability based on the betting scheme,  
a conditional probability assessment  $p(A|H)=x$ means that, for every real number $s$,  you are willing to pay 
an amount $s\cdot x$ and to receive $s$, or 0, or $s\cdot x$ (money back), according
to whether $AH$ is true, or $\widebar{A}H$ is true, or $\widebar{H}$
is true (bet called off), respectively. 
The random gain, which is  the difference between the (random) amount that you receive and the amount that you pay, is 
$G=(sAH+0\widebar{A}H+sx\widebar{H})-sx=
sAH+sx(1-H)-sx=
sH(A-x)$.

Given a probability function $p$ defined on an arbitrary family $\mathcal{K}$ of 
conditional events, consider a finite subfamily $\F = \{A_1|H_1, \ldots,A_n|H_n\} \subseteq \mathcal{K}$ and the vector
$\P=(x_1,\ldots, x_n)$, where $x_i = p(A_i|H_i)$ is the
assessed probability for the conditional event  $A_i|H_i$, $i=1,\ldots,n$.
With the pair $(\F,\P)$ we associate the random gain $G =
\sum_{i=1}^ns_iH_i(A_i - x_i)$. We denote by $\G_{\mathcal{H}_n}$ the set of values of $G$ restricted to $\H_n= H_1 \vee \cdots \vee H_n$, i.e., the set of values of $G$ when $\H_n$ is true. Then, we
recall below the notion of coherence in the context of the  {\em  betting scheme}.
\begin{definition}\label{COER-EV}{\rm
The function $p$ defined on $\mathcal{K}$ is coherent if and only if, $\forall n
		\geq 1$, $\forall \, s_1, \ldots,
		s_n$, $\forall \, \F=\{A_1|H_1, \ldots,A_n|H_n\} \subseteq \mathcal{K}$, it holds that: $min \; \G_{\mathcal{H}_n} \; \leq 0 \leq max \;
		\G_{\mathcal{H}_n}$. }
\end{definition}
In betting terms, the coherence of conditional probability assessments means that 
 in any finite combination of $n$ bets, 
after discarding the case where  all the bets are called off, the values of the random gain are  neither all positive nor all negative (i.e., \emph{no Dutch Book}).
 In particular,  coherence of $x=p(A|H)$ is defined by the condition $min \; \G_{H} \; \leq 0 \leq max \;
\G_{H}$, $\forall \, s$, where $\G_{H}$ is the set of values of $G$ restricted to $H$ (that is when the bet is not called off). Depending on the logical relations between $A$ and $H$ (with $H\neq \emptyset$),
 the set $\Pi$ of all coherent conditional probability assessments $x=p(A|H)$ is: 
 \begin{equation}\label{EQ:PI}
	\Pi= \left\{\begin{array}{ll}
		[0,1], &\mbox{if $\emptyset \neq AH \neq H$,}\\
		\{0\}, &\mbox{if $AH=\emptyset$ ,}\\
		\{1\}, &\mbox{if $AH=H$. }\\
	\end{array}
	\right.
\end{equation}  
 In numerical terms, once  $x=p(A|H)$ is assessed by the betting scheme, the indicator of $A|H$, denoted by the same symbol, is defined as  $1$, or $0$, or $x$, 
according
to whether $AH$ is true, or $\widebar{A}H$ is true, or $\widebar{H}$
is true. 
Then, by setting  $p(A|H)=x$, 
\begin{equation}\label{EQ:AgH}
	A|H=AH+x \no{H}=\left\{\begin{array}{ll}
		1, &\mbox{if $AH$ is true,}\\
		0, &\mbox{if $\no{A}H$ is true,}\\
		x, &\mbox{if $\no{H}$ is true.}\\
	\end{array}
	\right.
\end{equation}
Note that since the three-valued numerical entity  $A|H$ is defined by the betting scheme once the value $x=p(A|H)$ is assessed, the definition of (the indicator of) $A|H$ is not circular. The third  value of the random quantity  $A|H$  (subjectively) depends on the assessed probability  $p(A|H)=x$. 
Moreover, the value $x$ coincides with the corresponding conditional prevision, denoted by $\prev(A|H)$,   because
$\prev(A|H)=\prev(AH+x\no{H})=p(AH)+xp(\no{H})=p(A|H)p(H)+xp(\no{H})=xp(H)+xp(\no{H})=x$.

In the special case where  $AH=H$,  it follows by (\ref{EQ:PI}) that $x=1$ is the only coherent assessment for $p(A|H)$; then,
for the indicator $A|H$ it holds that 
\begin{equation}\label{EQ:AgH=1}
A|H=AH+x\no{H}=H+\no{H}=1, \;\; \mbox{ if }
AH=H. 
\end{equation}
In particular  (\ref{EQ:AgH=1}) holds when $A=\Omega$ (i.e., the sure event),  since   $\Omega\wedge H=H$ and hence
$
\Omega|H=H|H=1.
$
Likewise, if $AH=\emptyset$, 
it follows by (\ref{EQ:PI}) that $x=0$ is the only coherent assessment for $p(A|H)$; then,
\begin{equation}\label{EQ:AgH=0}
A|H=0+0\no{H}=0, \;\; \mbox{ if }
AH=\emptyset. 
\end{equation}
In particular  (\ref{EQ:AgH=0})  holds when $A=\emptyset$, since   $\emptyset\wedge H=\emptyset$ and hence
$
\emptyset|H=0.
$
We observe that conditionally on $H$ be true, 
for the (indicator of the) negation
it holds that $\overline{A|H}=\no{A}=1-A=
1-A|H.$ Conditionally on $H$ be false, by coherence, it holds that  $\overline{A|H}=p(\no{A}|H)=1-p(A|H)=1-A|H$. Thus, in all cases it holds that 
\begin{equation}\label{EQ:NEG}
\overline{A|H}=\no{A}|H=(1-A)|H=
1-A|H.
\end{equation}
We denote by $X$ a \emph{random quantity}, with a finite set of possible values. Given any event $H\neq \emptyset$, 
agreeing to the betting metaphor, if you  assess  the prevision  $\prev(X|H)=\mu$ means that for any given  real number $s$ you are willing to pay an amount $s \mu$ and to receive  $sX$, or $s\mu$, according  to whether $H$ is true, or  false (bet  called off), respectively. In particular, when $X$ is (the indicator of) an event $A$, then $\prev(X|H)=P(A|H)$.  The notion of coherence  can be generalized to the case of  prevision assessments on a family of conditional random quantities (see, e.g., \cite{GiSa20,SGOP20}). 
Given a random quantity $X$ and an event $H \neq \emptyset$,  
with prevision $\prev(X|H) = \mu$, likewise formula (\ref{EQ:AgH}) for the indicator of a conditional event, 
an extended notion of a conditional random quantity, denoted by the same symbol 
$X|H$, is defined as follows $X|H=XH+\mu\no{H}.$
We recall now  the notion of conjunction of two (or more) conditional events within the framework of conditional  random quantities in the setting of coherence
(\cite{GiSa13a,GiSa14,GiSa19,GiSa21IJAR}, for alternative approaches see also, e.g., \cite{Kaufmann2009,mcgee89}). 
Given a coherent probability assessment $(x,y)$ on $\{A|H,B|K\}$, we consider the random quantity $AHBK+x\no{H}BK+y\no{K}AH$ and we set $\prev[(AHBK+x\no{H}BK+y\no{K}AH)|(H\vee K)]=z$. Then  we define the  conjunction $(A|H)\wedge(B|K)$ as follows:
\begin{definition}\label{CONJUNCTION}{\rm Given a coherent prevision assessment 
		$p(A|H)=x$, $p(B|K)=y$, and $\prev[(AHBK+x\no{H}BK+y\no{K}AH)|(H\vee K)]=z$, the conjunction
		$(A|H)\wedge(B|K)$ is the conditional random quantity  defined as
		\begin{equation}\label{EQ:CONJUNCTION}
		\begin{array}{ll}
		(A|H)\wedge(B|K)=(AHBK+x\no{H}BK+y\no{K}AH)|(H\vee K) =\\
		
		=\left\{\begin{array}{ll}
		1, &\mbox{if $AHBK$ is true,}\\
		0, &\mbox{if $\no{A}H\vee \no{B}K$ is true,}\\
		x, &\mbox{if $\no{H}BK$ is true,}\\
		y, &\mbox{if $AH\no{K}$ is true,}\\
		z, &\mbox{if $\no{H}\no{K}$ is true}.
		\end{array}
		\right.
		\end{array}
		\end{equation}
}\end{definition}
Of course, $\prev[(A|H)\wedge(B|K)]=z$. 
Coherence requires that the Fréchet-Hoeffding bounds for prevision of the conjunction are preserved (\cite{GiSa14}), i.e, $\max\{x+y-1,0\}\leq z\leq \min\{x,y\}$, like in the case of unconditional events. Other  preserved properties are listed in \cite{GiSa21}.
We notice that if  conjunctions  of conditional events are defined as suitable conditional events (see, e.g., \cite{adams75,Cala17,CiDu12,CiDu13,GoNW91}),
 classical probabilistic properties are not preserved. In particular, the lower and upper probability bounds for the conjunction do not coincide with the above mentioned Fréchet-Hoeffding bounds (\cite{SUM2018S}).
Here, 
differently from conditional events which are three-valued objects, the conjunction $(A|H) \wedge (B|K)$
is not any longer  three-valued object, but  a five-valued object with values in $[0,1]$. We observe that $(A|H)\wedge (A|H)=A|H$ and $(A|H)\wedge (B|K)=(B|K)\wedge (A|H)$.
Moreover, if $H=K$, 
then  $(A|H)\wedge (B|H)=AB|H$.\\
For comparison with other approaches, like \cite{FlGH20}, see \cite[Section 9]{GiSa20}.

In analogy to  formula (\ref{EQ:AgH}), where the indicator of  a conditional event ``$A$ given $H$'' is defined as $A|H=A\wedge H+p(A|H)\no{H}$, 
the iterated conditional  ``$B|K$ given $A|H$'' is defined as follows (see, e.g., \cite{GiSa13c,GiSa13a,GiSa14}):
\begin{definition}[Iterated conditioning]
	\label{DEF:ITER-COND} Given any pair of conditional events $A|H$ and $B|K$, with $AH\neq \emptyset$, the iterated
	conditional $(B|K)|(A|H)$ is defined as the conditional random quantity
	$(B|K)|(A|H) = (A|H) \wedge (B|K) + \mu \no{A}|H$,
	where
	$\mu =\mathbb{P}[(B|K)|(A|H)]$.
\end{definition}
	Notice that we assume  $AH\neq \emptyset$ to avoid trivial cases of iterated conditionals. 
Specifically, similar as in the  three-valued notion of a conditional event $A|H$ where the antecedent $H$ must not be impossible (i.e.,  $H$ must not coincide with the constant 0),  the iterated conditional $(B|K)|(A|H)$ requires that the antecedent $A|H$ must  not be constant and equal to 0 (this happens when $AH\neq \emptyset$). Furthermore, we recall that the compound prevision theorem is preserved, that is 
$\prev[(A|H) \wedge (B|K)]=\prev[(B|K)|(A|H)]p(A|H)$.
\section{Approach 1: Connexive Principles and Default Reasoning}
In order to validate the connexive principles we interpret
a conditional $A\rightarrow C$ by the 
default $A\normally C$, where $A$ and $C$ are two events (with $A\neq \bot$).
A default $A\normally C$ can be read as \emph{$C$ is a plausible consequence of $A$} and  is interpreted by 
the probability constraint $p(C|A)=1$
(\cite{gilio16}).\footnote{
According to $\varepsilon$-semantics 
(see, e.g., \cite{adams75,pearl89}) a default $A\normally C$ is interpreted  by $p(C|A)\geq 1-\varepsilon$, with $\varepsilon>0$ and $p(A)>0$. 
Gilio introduced a coherence-based probability semantics for defaults by also allowing $\varepsilon$ and $p(A)$ to be zero (\cite{gilio02}). In this context,  defaults in terms of probability 1 can be used to give a alternative definition of p-entailment which  preserve  the usual non-monotonic inference rules like those  of System P (\cite{biazzo02,gilio02,Gilio2013,gilio13}, see also \cite{coletti02,Coletti2015}). For the psychological plausibility of the coherence-based semantics of non-monotonic reasoning, see, e.g., \cite{pfeifer,pfeifer09b,pfeifer10a,pfeifertulkki17}.} The conditional  $A \rightarrow  \lsim C$ is interpreted by the default $A\normally \no{C}$. Likewise, $\lsim A \rightarrow  C$ is interpreted by $\no{A}\normally C$.
A negated conditional  $\lsim (A\rightarrow C)$ is interpreted by the negated default $\lsim (A \normally C)$ (\emph{it is not the case that: $C$ is a plausible consequence of $A$};
also  denoted by
$A\nnormally C$ in 
\cite{freund1991,kraus90}) that is  $p(C|A)\neq 1$, which corresponds to the wide scope negation of negating conditionals (\cite{gilio16}).

The conjunction of two conditionals, denoted by
$(A\rightarrow B)\wedge (C\rightarrow D)$, is interpreted by the sequence of their associated  defaults $(A\normally B,C\normally D)$, which  represents in probabilistic terms  the constraint
$(p(B|A)=1,p(D|C)=1)$, that is
$(p(B|A),p(D|C))=(1,1)$.
Then, the negation of the conjunction of two conditionals, denoted by $\lsim ((A\rightarrow B)\wedge (C\rightarrow D))$ (i.e.,  in terms of defaults $\lsim (A\normally B)\wedge (C\normally  D)$) 
is interpreted by the
negation of the 
probabilistic constraint $(p(B|A),p(D|C))=(1,1)$, that is  $(p(B|A),p(D|C))\neq (1,1)$. Table~\ref{TAB:BS}  summarizes the  interpretations.

We now introduce the definition of validity for non-iterated connexive principles (e.g., (AT), (AT$'$), (AB)).
\begin{definition}\label{DEF:V1}
We say that a non-iterated connexive principle is valid if and only if the probabilistic constraint associated with  the connexive principle is satisfied by  every  coherent assessment on the involved conditional events. 
\end{definition}
In the next paragraphs we check the validity in terms of  Definition \ref{DEF:V1} of the non-iterated connexive principles in Table \ref{TAB:CP}. 
\paragraph{Aristotle's Thesis (AT): $\lsim (\lsim A \rightarrow A)$.}
We assume that $\no{A}\neq \bot$ and we interpret the principle $\lsim (\lsim A \rightarrow A)$ by the negated default $\lsim (\no{A}\normally A)$ with the following associated probabilistic constraint: $p(A|\no{A})\neq 1$.
We observe that $p(\no{A}|A)=0$ is the unique precise coherent assessment on $\no{A}|A$. Then, 
(AT) is valid because every coherent precise assessment $p(A|\no{A})$ is such that $p(A|\no{A})\neq 1$.\\[-2em]  
\paragraph{Aristotle's Thesis $'$ (AT $'$): $\lsim (A \rightarrow \lsim A)$.}
 Like (AT),   (AT)$'$ can be validated.\\[-2em] 
\paragraph{Abelard's Thesis (AB): $\lsim((A \rightarrow B) \wedge (A \rightarrow\lsim B))$.} We assume that $A\neq \bot$.
The structure of this principle  is formalized by $\lsim((A \normally B) \wedge (A \normally  \no{B}))$ which expresses  the constraint $(p(B|A),p(\no{B}|A))\neq (1,1)$. 
We recall that  coherence requires  $p(B|A)+p(\no{B}|A)=1$. Then, (AT) is valid because each coherent assessment on $(B|A,\no{B}|A)$ is necessarily of the form $(x,1-x)$, with $x\in[0,1]$, which of course satisfies   $(p(B|A),p(\no{B}|A))\neq (1,1)$.\\[-2em]  \paragraph{Aristotle's Second Thesis (AS): $\lsim((A \rightarrow B) \wedge (\lsim A \rightarrow B))$.} 
We assume that $A \neq \emptyset$ and $\no{A} \neq \emptyset$.
The structure of this principle  is formalized by
$\lsim((A \normally B) \wedge (\no{A} \normally  B))$ which expresses  the constraint $(p(B|A),p(B|\no{A}))\neq (1,1)$. We recall that, given two logically independent events $A$ and $B$ every assessment $(x,y)\in[0,1]^2$ on $(B|A,B|\no{A})$ is coherent. In particular, $(p(B|A),p(B|\no{A}))= (1,1)$ is a coherent assessment which does not satisfy the probabilistic constraint $(p(B|A),p(B|\no{A}))\neq (1,1)$. Thus, (AS) is not valid.

Concerning iterated connexive principles (e.g. (BT), (BT$'$)), we interpret the main connective ($\rightarrow$) as the implication ($\Rightarrow$)  from the  probabilistic constraint on the antecedent to the probabilistic constraint on the conclusion. 
Then, for instance, the iterated conditional
$(A\rightarrow B)\rightarrow (C\rightarrow D)$ is interpreted by the implication  $A\normally B \Rightarrow C\normally D$, that is
$p(B|A)=1  \Rightarrow  p(D|C)=1$.
We now define validity for iterated connexive principles.
\begin{definition}\label{DEF:V1ITER}
An iterated connexive principle $\bigcirc  \Rightarrow \Box$ is valid  if and only if the probabilistic constraint of the conclusion $\Box$ is satisfied by  every  coherent extension to the conclusion from any coherent probability assessment satisfying the constraint of the premise $\bigcirc$. 
\end{definition}
\begin{table}[!ht]	
\vspace{-1.3em}	
\begin{center}
	\begin{tabular}{ccc}	
	\hline 
		Conditional object  & Default  & Probabilistic interpretation  \\ \hline	
		$A \rightarrow C$ & $A\normally C$ & $p(C|A)=1$ \\
       $\lsim (A \rightarrow C)$ & $\lsim (A\normally C)$ & $p(C|A)\neq 1$\\   
       $(A\rightarrow B)\wedge (C\rightarrow D)$ &\quad $(A\normally B, C\normally D)$\quad & $(p(B|A),p(D|C))=(1,1)$\\
      $\lsim((A\rightarrow B)\wedge (C\rightarrow D))$ &\quad $\lsim(A\normally B, C\normally D)$\quad & $(p(B|A),p(D|C))\neq (1,1)$\\
      $(A\rightarrow B)\rightarrow (C\rightarrow D)$ & $A\normally  B \Rightarrow C\normally  D$ & $p(B|A)=1\, \Rightarrow\, p(D|C)=1$\\ 
\hline 
	\end{tabular}
\end{center}
\vspace{-0.5em}
\caption{Probabilistic interpretations of  logical operation on conditionals  in terms of defaults or negated defaults.}
\label{TAB:BS} %
\vspace{-0.5em}
\end{table}	\vspace{-2.3em}
We check the validity in terms of Definition \ref{DEF:V1ITER} of the iterated connexive principles in Table \ref{TAB:CP}.\\[-2em]  
\paragraph{Boethius' Thesis (BT): $(A \rightarrow B) \rightarrow \lsim (A
  \rightarrow \lsim B)$.}
We assume that $A \neq \emptyset$.
This
 is interpreted by the implication
 $A\normally  B \Rightarrow  \lsim (A\normally  \no{B})$, that is   $p(B|A)=1\, \Rightarrow\, p(\no{B}|A)\neq 1$.
 We observe that, by setting $p(B|A)=1$, 
 $p(\no{B}|A)=1-p(B|A)=0$ is the unique coherent extension  to $\no{B}|A$. Then, as $p(B|A)=1 \Rightarrow p(\no{B}|A)=0\neq 1$, the iterated connexive principle (BT) is valid.\\[-2em]  
 \paragraph{Boethius' Thesis $'$ (BT $'$): $(A \rightarrow \lsim B) \rightarrow
   \lsim(A \rightarrow B)$.} Like (BT), it can be shown that (BT$'$) is valid too.\\[-2em] 
  \paragraph{Reversed Boethius' Thesis (RBT): $\lsim (A
  \rightarrow \lsim B) \rightarrow(A \rightarrow B)$.}
  We assume that $A\neq \bot$.
This
 is interpreted by
 $\lsim (A\normally  \no{B}) \Rightarrow   A\normally  B$, that is   $p(\no{B}|A)\neq1\, \Rightarrow\, p(B|A)= 1$.
 We observe that, 
 by setting $p(\no{B}|A)=x$ it holds that  $p(B|A)=1-x$ is the unique coherent extension  to $B|A$. In particular by choosing $x\in]0,1[$, it holds that $p(\no{B}|A)\neq 1$ and $p(B|A)\neq 1$.
 Thus,  $p(\no{B}|A)\neq 1 \nRightarrow p(B|A)= 1$ and hence  (RBT) is not valid.\\[-2em]  
 \paragraph{Reversed Boethius' Thesis $'$ (RBT $'$): $\lsim (A
  \rightarrow  B) \rightarrow(A \rightarrow \lsim B)$.}
 Like (RBT), it can be shown that (RBT$'$) is not valid too.\\[-2em] 
\paragraph{Boethius Variation (B3): $(A \rightarrow B) \rightarrow  \lsim (\lsim A
  \rightarrow  B)$.}
 We assume that $A\neq \bot$ and $\no{A}\neq \bot$. 
 This
 is interpreted by
 $A\normally  B \Rightarrow   \lsim(\no{A}\normally  B)$, that is   $p(B|A)=1\, \Rightarrow\, p(B|\no{A})\neq  1$.
 We observe that, by setting $p(B|A)=1$, any value $p(B|\no{A})\in[0,1]$ is a coherent extension  to $B|\no{A}$, because the assessment $(1,y)$ on $(B|A,B|\no{A})$ is coherent for every  $y\in[0,1]$. In particular the assessment $(1,1)$ on $(B|A,B|\no{A})$ is coherent. Therefore, as $p(B|A)=1 \nRightarrow p(B|\no{A})\neq 1$, 
 (B3) is not valid.\\[-2em]  
\paragraph{Boethius Variation (B4): $(\sim A \rightarrow B) \rightarrow  \lsim (A
  \rightarrow  B)$.} 
Like (B3), it can be shown that (B4) is not valid too. 

We summarize the results of this section in Table \ref{TAB:CC1}.
\begin{table}[!h]
	\begin{center}
\begin{scriptsize}
		\begin{tabular}{ccccc}	
			\hline 
			Name
			&	
			Connexive principle   & Default  &Probabilistic constraint & Validity   \\\hline	&&\\[-.8em] 
			(AT) &$\lsim (\lsim A \rightarrow A)$ & $\lsim(\no{A} \normally A)$ & $p(A|\no{A})\neq 1$ & yes\\	
			
			&&\\[-.8em] 
			
			(AT') &$\lsim ( A \rightarrow \lsim A)$ & 
			$\lsim (A \normally \no{A})$&
			$p(\no{A}|A)\neq 1$& yes\\
			
			&&\\[-.8em] 
			
			(AB) &$\lsim((A \rightarrow B) \wedge (A \rightarrow\lsim B))$& $ \lsim(A \normally B, A \normally \no{B})$ & $(p(B|A),p(\no{B}|A))\neq (1,1)$& yes\\
			
			&&\\[-.8em] 
			
			(AS) & $\lsim((A \rightarrow B) \wedge (\lsim A \rightarrow B))$&
			$ \lsim(A \normally B, \no{A} \normally B)$&
		$(p(B|A),p(B|\no{A}))\neq (1,1)$& no\\
			
			&&\\[-.8em] 
			
			(BT) & $(A \rightarrow B) \rightarrow \lsim (A
			\rightarrow \lsim B)$&
			$A\normally B \Rightarrow \lsim (A \normally \no{B})$
			&
			$p(B|A)=1 \Rightarrow p(\no{B}|A)\neq 1$& yes\\
			
			&&\\[-.8em] 
			
			(BT') & $(A \rightarrow \lsim B) \rightarrow
			\lsim(A \rightarrow B)$ & \quad
			$A\normally \no{B} \Rightarrow \lsim (A \normally B)$
			\quad&
			$p(\no{B}|A)=1 \Rightarrow p(B|A)\neq 1$& yes\\
			
			&&\\[-.8em] 
			
			(RBT) & $\lsim (A
			\rightarrow \lsim B) \rightarrow(A \rightarrow B)$ &
			$\lsim (A \normally \no{B})\Rightarrow  A\normally B  $
			&
			$p(\no{B}|A)\neq 1 \Rightarrow p(B|A)= 1$& no\\
			
			&&\\[-.8em] 
			
			(RBT') & $ \lsim(A \rightarrow B)\rightarrow(A \rightarrow \lsim B)$&
			$\lsim(A\normally B) \Rightarrow  A \normally \no{B}$
			&
			$p(B|A)\neq 1 \Rightarrow p(\no{B}|A)= 1$& no\\
			
			&&\\[-.8em] 
			
			(B3) & $(A \rightarrow B) \rightarrow  \lsim (\lsim A
			\rightarrow  B)$ &
			$A\normally B \Rightarrow \lsim( \no{A} \normally B)$
			&
			$p(B|A)= 1 \Rightarrow p(B|\no{A})\neq 1$& no\\
			
			&&\\[-.8em] 
			
			(B4)& $(\sim A \rightarrow B) \rightarrow  \lsim (A
			\rightarrow  B)$&
			$\no{A}\normally B \Rightarrow \lsim( A \normally B)$&
			$p(B|\no{A})= 1 \Rightarrow p(B|A)\neq  1$& no\\
			
			&&\\[-.8em] 
			
			\hline 
		\end{tabular}
\end{scriptsize}		
	\end{center}
	\vspace{-0.5em}
	\caption{Connexive principles in the framework of defaults and probabilistic constraints (Approach 1).}
	\label{TAB:CC1}
	\vspace{-0.5em}%
\end{table}	
\section{Approach 2: Connexive Principles and Compounds of Conditionals}
In this section 
we analyze connexive principles  within the theory of logical operations among conditional events. 
Specifically, we analyze connections between antecedents and consequents in terms of constraints on compounds of conditionals and iterated conditionals.
In this second approach, a basic conditional $A\rightarrow C$ is interpreted as a conditional event $C|A$ (instead of a  probabilistic constraint on conditional events) which is a three-valued object: $C|A\in\{1,0,x\}$, where $x=p(C|A)$. 
The negation $\lsim (A\rightarrow C)$ is interpreted by $\no{C}|A$ (which is the narrow scope negation of negating conditionals). Then, $\lsim (A\rightarrow \lsim C)$ amounts to $\overline{\no{C}|A}$ which coincides with $C|A$. 
We recall that logical operations among conditional events do not yield  a conditional event, rather they yield  conditional random quantities with more than three possible values (see, e.g.,  \cite{GiSa14}). Then, we interpret the results of the logical operations in the connexive principles  by suitable conditional random quantities. 
In particular, the conjunction  $(A\rightarrow B)\wedge (C\rightarrow D)$ (resp., $\lsim ((A\rightarrow B)\wedge (C\rightarrow D))$) is  interpreted by $(B|A)\wedge (D|C)$ (resp., by $\overline{(B|A)\wedge (D|C)}$), and the iterated conditional $(A\rightarrow B)\rightarrow (C\rightarrow D)$ is interpreted by $(D|C)|(B|A)$.
Moreover, we define validity of connexive principles within Approach~2. 
\begin{definition}\label{DEF:V2}
 A  connexive principle is valid  if and only if the associated conditional random quantity is constant and equal to 1.
\end{definition}
We now check the validity of the connexive principles in Table \ref{TAB:CP} according to Definition \ref{DEF:V2}.

\paragraph{Aristotle's Thesis (AT): $\lsim (\lsim A \rightarrow A)$.}
We interpret the principle $\lsim (\lsim A \rightarrow A)$ by the negation of the  conditional event $A|\no{A}$, that is by $\overline{A|\no{A}}$, where $\no{A}\neq \emptyset$.
 Then, based on equations (\ref{EQ:NEG}) and  (\ref{EQ:AgH=1}), it follows that
$
\overline{A|\no{A}}=1-A|\no{A}=\no{A}|\no{A}=1
$.
Therefore,  (AT) is valid because the conditional random quantity $\overline{A|\no{A}}$, which also coincides with the conditional event $\no{A}|\no{A}$,  is constant and equal to 1.\\[-2em] 
\paragraph{Aristotle's Thesis $'$ (AT $'$): $\lsim ( A \rightarrow \lsim A)$.}
We interpret the principle $\lsim (A \rightarrow \lsim A)$ by the negation of the  conditional event $\no{A}|A$, that is by $\overline{\no{A}|A}$, where $A\neq \emptyset$. 
Like in (AT), it holds that  
$
\overline{\no{A}|A}=1-\no{A}|A=A|A=1
$,
which validates (AT $'$). Notice that, (AT $'$) also follows from (AT) when $A$ is replaced by $\no{A}$ (of course $\no{\no{A}}=A$). \\[-2em]

\paragraph{Abelard's Thesis (AB): $\lsim((A \rightarrow B) \wedge (A \rightarrow\lsim B))$.}
The structure of this principle  is formalized by the conditional random quantity $\overline{(B|A)\wedge (\no{B}|A)}$, where $A\neq \emptyset$. 
We observe that $(B|A)\wedge (\no{B}|A)=(B\wedge \no{B})|A=\emptyset|A$. Then, 
$\overline{(B|A)\wedge (\no{B}|A)}=\overline{\emptyset|A}=
\no{\emptyset}|A=\Omega|A=1.$     
Therefore,  (AB) is valid.\\[-2em]

\paragraph{Aristotle's Second Thesis (AS): $\lsim((A \rightarrow B) \wedge (\lsim A \rightarrow B))$.}
The structure of this principle  is formalized by the  random quantity $\overline{(B|A) \wedge (B|\no{A})}$, where $A \neq \emptyset$ and $\no{A} \neq \emptyset$.
By setting $p(B|A)=x$ and $p(B|\no{A})=y$, it follows that  \cite{GiSa13a,GiSa19b}
\[
(B|A) \wedge (B|\no{A})=(B|A) \cdot  (B|\no{A})
=\left\{
\begin{array}{ll}
0,& \mbox{ if } A\no{B}  \vee   \no{A}\no{B}\mbox{ is  true}, \\
y,     & \mbox{ if } AB \mbox{ is  true}, \\
x,               & \mbox{ if } \no{A}B \mbox{ is  true}.
\end{array}
\right.
\]
Then, $\overline{(B|A) \wedge (B|\no{A})}=$$1-(B|A) \wedge (B|\no{A})=1-(y AB+x\no{A}B)$,
which is not constant and can therefore  not necessarily be equal to  1. In particular, by choosing the  coherent assessment $x=y=1$, it follows that $\overline{(B|A) \wedge (B|\no{A})}=1-AB-\no{A}B=1-B=\no{B}$, which is not necessarily equal to 1 as it could be either 1 or 0, according to whether $\no{B}$ is  true or false, respectively. Therefore, (AS) is not valid.
Moreover,  by setting $\prev[(B|A) \wedge (B|\no{A})]=\mu$, it holds that
$
\mu=y\,p(AB)+x\,p(\no{A}B)=$$y\,p(B|A)p(A)+x\,p(B|\no{A})p(\no{A})=$$xy\,p(A)+xy\,p(\no{A})=xy$. 
Then,
$\prev[\overline{(B|A) \wedge (B|\no{A})}]=1-xy$. We also observe that, in the special case where $x=y=0$, it follows that $\overline{(B|A) \wedge (B|\no{A})}=1$.\\[-2em]

\paragraph{Boethius' Thesis (BT): $(A \rightarrow B) \rightarrow \lsim (A
  \rightarrow \lsim B)$.}
This principle 
 is formalized by the iterated conditional  $(\overline{\no{B}|A})|(B|A)$, with $AB\neq \emptyset$. We recall that  $(\overline{\no{B}|A})=B|A$. Then  $(\overline{\no{B}|A})|(B|A)=(B|A)|(B|A)$. Moreover, by setting $p(B|A)=x$ and  $\prev[(B|A)|(B|A)]=\mu$, it holds that 
 \[
 \begin{array}{ll}
 (B|A)|(B|A)=(B|A)\wedge (B|A)+\mu(1-B|A)=(B|A)+\mu(1-B|A)=\\
 =\left\{
\begin{array}{ll}
1,& \mbox{ if } AB \mbox{ is  true}, \\
\mu,     & \mbox{ if } A\no{B}, \mbox{ is  true}, \\
x+\mu(1-x),              & \mbox{ if } \no{A} \mbox{ is  true}.
\end{array}
\right.
\end{array}
 \]
 By linearity of prevision it holds that $\mu=x+\mu(1-x)$. Then,
 \[
 (B|A)|(B|A)=(B|A)+\mu(1-B|A)\left\{
\begin{array}{ll}
1,& \mbox{ if } AB \mbox{ is  true}, \\
\mu,     & \mbox{ if } \no{A} \vee \no{B} \mbox{ is  true}. \\
\end{array}
\right.
 \]
 Then, by coherence it must be that $\mu=1$ and hence  (\cite[Remark 2]{GiSa14}, see also \cite[Section 3.2]{GiSa21})
\begin{equation}\label{EQ:BgAgBgA}
(B|A)|(B|A)=1.
\end{equation} 
 Therefore $(\overline{\no{B}|A})|(B|A)$ is constant and equal to 1 and hence   (BT) is valid.\\[-2em] 

\paragraph{Boethius' Thesis $'$ (BT $'$): $(A \rightarrow \lsim B) \rightarrow
  \lsim(A \rightarrow B)$.}
This principle is formalized by the iterated conditional $(\overline{B|A})|(\no{B}|A)$, where $A\no{B}\neq \emptyset$. By observing that $\overline{B|A}=\no{B}|A$, it follows that $(\overline{B|A})|(\no{B}|A)=(\no{B}|A)|(\no{B}|A)$ which is constant and equal to 1 because of (\ref{EQ:BgAgBgA}). Therefore, (BT$'$) is valid.\\[-2em]
\paragraph{Reversed Boethius' Thesis (RBT): $\lsim (A
  \rightarrow \lsim B) \rightarrow(A \rightarrow B)$.}
This principle is formalized by the iterated conditional $(B|A)|(\overline{\no{B}|A})$, where $AB\neq \emptyset$. As $(\overline{\no{B}|A})=B|A$,  it follows 
from (\ref{EQ:BgAgBgA})
that $(B|A)|(\overline{\no{B}|A})=(B|A)|(B|A)=1$. Therefore, (RBT) is valid.\\[-2em]
\paragraph{Reversed Boethius' Thesis $'$ (RBT $'$): $ \lsim(A \rightarrow B)\rightarrow(A \rightarrow \lsim B)$.}
This principle is formalized by the iterated conditional $(\no{B}|A)|(\overline{B|A})$, where $A\no{B}\neq \emptyset$. As $(\overline{B|A})=\no{B}|A$,  it follows from (\ref{EQ:BgAgBgA}) that $(\no{B}|A)|(\overline{B|A})=(\no{B}|A)|(\no{B}|A)=1$. Therefore (RBT$'$) is validated.\\[-2em]

\paragraph{Boethius Variation (B3): $(A \rightarrow B) \rightarrow  \lsim (\lsim A
  \rightarrow  B)$.} This principle is formalized by the iterated conditional $(\overline{B|\no{A}})|(B|A)$, where $AB\neq \emptyset$.
  We observe that $(\overline{B|\no{A}})|(B|A)=(\no{B}|\no{A})|(B|A)$, because $\overline{B|\no{A}}
  =\no{B}|\no{A}$.
By setting $p(B|A)=x$, $p(\no{B}|\no{A})=y$, and  $\prev[(\no{B}|\no{A})|(B|A)]=\mu$, it holds that 
 \[
 (\no{B}|\no{A})|(B|A)=(\no{B}|\no{A})\wedge (B|A)+\mu(1-B|A)=
 \left\{
\begin{array}{ll}
y,& \mbox{ if } AB \mbox{ is  true}, \\
\mu,     & \mbox{ if } A\no{B} \mbox{ is  true}, \\
\mu(1-x),              & \mbox{ if } \no{A}B \mbox{ is  true},\\
x+\mu(1-x),              & \mbox{ if } \no{A}\no{B} \mbox{ is  true},
\end{array}
\right.
 \]
which is not constant and can therefore  not necessarily be equal to  1.
For example, if we  choose the  coherent assessment $x=y=1$, it follows that  \[
 (\no{B}|\no{A})|(B|A)=(\no{B}|\no{A})\wedge (B|A)+\mu(1-B|A)=
 \left\{
\begin{array}{ll}
1,& \mbox{ if } AB \mbox{ is  true}, \\
\mu,     & \mbox{ if } A\no{B} \mbox{ is  true}, \\
0,              & \mbox{ if } \no{A}B \mbox{ is  true},\\
1,              & \mbox{ if } \no{A}\no{B} \mbox{ is  true},
\end{array}
\right.
 \]
which is not constant and equal to 1. Therefore, (B3) is not valid.\\[-2em]

\paragraph{Boethius Variation (B4): $(\sim A \rightarrow B) \rightarrow  \lsim (A
  \rightarrow  B)$.}    This principle is formalized by the iterated conditional $(\overline{B|A})|(B|\no{A})$, where $\no{A}B\neq \emptyset$.
  We observe that $(\overline{B|A})|(B|\no{A})$ is not constant and not necessarily equal to 1 because it is equivalent to (B3) when $A$ is replaced by $\no{A}$. Therefore, (B4) is not valid.

Connexive principles and their interpretation in terms of compound or iterated conditionals are illustrated 
in Table \ref{TAB:CC}. 
\begin{table}[!ht]	
\begin{center}\begin{scriptsize}
	\begin{tabular}{ccrlc}	
	\hline 
Name
&	
		Connexive principle   & Interpretation & Value &Validity   \\\hline	&&\\[-.9em] 
(AT) &$\lsim (\lsim A \rightarrow A)$ & $\overline{ A|\no{A}} $&$=1$ & yes\\	

&&\\[-1em] 

(AT$'$) &$\lsim ( A \rightarrow \lsim A)$ & $\overline{\no{A} | A}$&$=1$& yes\\

&&\\[-1em] 

(AB) &$\lsim((A \rightarrow B) \wedge (A \rightarrow\lsim B))$& $\overline{(B|A) \wedge (\no{B}|A)}$&$=1$& yes\\

&&\\[-1em] 

 (AS) & $\lsim((A \rightarrow B) \wedge (\lsim A \rightarrow B))$&$\overline{(B|A) \wedge (B|\no{A})}$&$\neq 1$& no\\

&&\\[-1em] 

 (BT) & $(A \rightarrow B) \rightarrow \lsim (A
  \rightarrow \lsim B)$&$\overline{(\no{B}|A)}|(B|A)$&$=1$& yes\\
 
 &&\\[-1em] 

 (BT$'$) & $(A \rightarrow \lsim B) \rightarrow
  \lsim(A \rightarrow B)$ &$\overline{(B|A)}|(\no{B}|A)$&$=1$& yes\\

&&\\[-1em] 

(RBT) & $\lsim (A
  \rightarrow \lsim B) \rightarrow(A \rightarrow B)$ &$(B|A)|\overline{(\no{B}|A)}$&$=1$& yes\\

&&\\[-1em] 

(RBT$'$) & $ \lsim(A \rightarrow B)\rightarrow(A \rightarrow \lsim B)$&$(\no{B}|A)|\overline{(B|A)}$&$=1$& yes\\

&&\\[-1em] 

  (B3) & $(A \rightarrow B) \rightarrow  \lsim (\lsim A
  \rightarrow  B)$ &$\overline{(B|\no{A})}|(B|A)$&$\neq 1$& no\\
  
  &&\\[-1em] 
  
(B4)& $(\sim A \rightarrow B) \rightarrow  \lsim (A
  \rightarrow  B)$&$\overline{(B|A)}|(B|\no{A})$&$\neq 1$& no\\
 
  &&\\[-1em] 

\hline 
	\end{tabular}\end{scriptsize}
\end{center}
\vspace{-0.5em}
\caption{Connexive principles in the framework of compounds of conditionals and iterated conditionals (Approach 2). Value denotes whether the conditional random quantity is constant and equal to 1.}
\label{TAB:CC} %
\vspace{-0.8em}
\end{table}	
\section{Concluding Remarks}\vspace{-1em}
We presented two approaches to investigate connexive principles. Connexivity is interpreted by  in terms of probabilistic constraints on conditional events  (in  the  sense  of  defaults,  or  negated  defaults) in  Approach~1. Within this approach we showed that the connexive principles (AT), (AT$'$), (AB), (BT), and (BT$'$) are valid, whereas (AS), (RBT), (RBT$'$), (B3), and (B4) are not valid (see Table~\ref{TAB:CC1}). In Approach~2 connexivity is interpreted by constraints on compounds of conditionals and iterated conditionals. Here, we demonstrated that, like in Approach~1,  (AT), (AT$'$), (AB), (BT), and (BT$'$) are valid, whereas (AS),  (B3), and (B4) are not valid. Contrary to Approach~1, (RBT) and (RBT$'$) are valid in Approach~2 (see Table~\ref{TAB:CC}). 

Approach~1 is characterized by employing concepts from coherence-based probability theory and probabilistic interpretations of defaults and negated defaults.
Conditionals, interpreted  as  defaults, are negated 
 by the wide scope negation. We gave two notions of validity, namely  for non-iterated and iterated connexive principles, respectively. 
Approach~2  allows for dealing with logical operations on conditional events 
and avoids (see, e.g., \cite{SGOP20}) the well known
Lewis' triviality results (see, e.g.,
\cite{lewis76}).
It therefore 
offers  a more unified approach to connexive principles, which is reflected by a unique  definition of validity for both, iterated and non-iterated  connexive principles. Moreover, Approach~2  negates conditionals by the narrow scope negation. 
Thus, validity depends on  how conditionals and negation are defined.

One might wonder why neither of the two approaches validates all connexive principles. Of course, we have shown by proofs why, for instance, (AS) is not valid in both approaches. Apart from the insight obtained from our proofs, this is not surprising since also not all connexive principles are valid in all systems of connexive logic. Moreover,  some rules which are valid in classical logic (e.g., transitivity, contraposition, and premise strengthening) are not valid in probability logic, while, for example,  the rules of the basic nonmonotonic System P are valid within coherence-based probability logic (\cite{coletti02,gilio02,gilio16,GiSa19,GiPS20}).

We have shown that coherence-based probability logic  offers a rich language to
investigate the validity of various connexive principles. Future work will be devoted to investigations on other intuitively plausible logical principles contained in  alternative and non-classical logics.  \\ \ 

\noindent
{\bf Acknowledgments.} Thanks to  four anonymous reviewers for useful comments. 
\bibliographystyle{splncs04} 
\bibliography{ecsqaru2021PS} 

\end{document}